\documentclass[12pt]{article}
\usepackage{amsmath,amssymb,amsthm,amsfonts}

\newtheorem{thm}{Theorem}[section]
\newtheorem{lemma}[thm]{Lemma}
\newtheorem{prop}[thm]{Proposition}

\newtheorem{rem}[thm]{Remark}

\newcommand{\R}{{\mathbb{R}}}

\newcommand{\Q}{{\mathbb{Q}}}
\newcommand{\T}{{\mathbb{T}}}
\newcommand{\Z}{{\mathbb{Z}}}
\newcommand{\N}{{\mathbb{N}}}
\newcommand{\C}{{\mathbb{C}}}
\newcommand{\SP}{{\mathbb{S}}}

\newcommand{\cA}{{\mathcal{A}}}
\newcommand{\cB}{{\mathcal{B}}}
\newcommand{\bcB}{{\bar{\mathcal{B}}}}

\newcommand{\cF}{{\mathcal{F}}}

\newcommand{\cK}{{\mathcal{K}}}
\newcommand{\tcK}{\widetilde{\mathcal{K}}}

\newcommand{\cR}{{\mathcal{R}}}
\newcommand{\cS}{{\mathcal{S}}}

\newcommand{\Ca}{{\hbox{\it Cal}}}

\newcommand{\tLambda}{{\widetilde{\Lambda}}}

\newcommand{\tomega}{{\tilde{\omega}}}

\newcommand{\tM}{{\widetilde{M}}}

\newcommand{\Ham}{{\hbox{\it Ham\,}}}

\newcommand{\tHam}{\widetilde{\hbox{\it Ham}\, }}
\newcommand{\Symp}{{\hbox{\it Symp} }}

\newcommand{\Qed}{\hfill \qedsymbol \medskip}

\begin{document}

\title{Symplectic Quasi-States
and Semi-Simplicity of Quantum Ho\-mo\-lo\-gy}

\renewcommand{\thefootnote}{\alph{footnote}}

\author{\textsc Michael Entov$^{a}$\ and Leonid
Polterovich$^{b}$ }

\footnotetext[1]{Partially supported by E. and J. Bishop Research
Fund and by the Israel Science Foundation grant $\#$ 881/06.}
\footnotetext[2]{Partially supported by the Israel Science
Foundation grant $\#$ 11/03.}

\date{}

\maketitle

\begin{abstract}
We review and streamline our previous results and the results of
Y.~Ostrover on the existence of Calabi quasi-morphisms and
symplectic quasi-states on symplectic manifolds with semi-simple
quantum homology. As an illustration, we discuss the case of
symplectic toric Fano 4-manifolds. We present also new results due
to D.~McDuff: she observed that for the existence of
quasi-morphisms/quasi-states it suffices to assume that the
quantum homology contains a field as a direct summand, and she
showed that this weaker condition holds true for one point
blow-ups of non-uniruled symplectic manifolds.
\end{abstract}



\renewcommand{\thefootnote}{\arabic{footnote}}

\section{Symplectic quasi-states and Calabi quasi-morphisms} \label{sec-introduction}

\medskip
\noindent{\sc Symplectic quasi-states:}  Let $(M^{2n},\omega)$ be a
closed connected $2n$-di\-men\-sio\-nal symplectic manifold.  A {\it
symplectic quasi-state} on $M$ is a (possibly non-linear!)
real-valued functional $\zeta$ on the space $C(M)$ of all continuous
functions on $M$ which  satisfies the following conditions:

\medskip
\noindent {\bf Quasi-linearity:} $\zeta (F+\lambda G)=\zeta
(F)+\lambda \zeta (G)$ for all $\lambda\in\R$ and for all
functions $F,G\in C^\infty (M)$ which commute with respect to the
Poisson brackets: $\{F,G\} = 0$.

\medskip
\noindent {\bf Monotonicity:} $\zeta (F) \leq \zeta (G)$ if  $F
\leq G$.

\medskip
\noindent {\bf Normalization:} $\zeta (1) = 1$.

\medskip
\noindent Non-linear quasi-states are known to exist on certain
symplectic manifolds. For instance, we will prove in this paper
that they exist on any symplectic toric Fano four-manifold
equipped with a rational symplectic form. Here, denoting by $H_2^S
(M)$ the group of spherical integral homology classes of $M$, we
will say that $\omega$ is {\it rational} if the subgroup $[\omega]
(H_2^S (M))\subset \R$ is discrete. The list of symplectic toric
Fano four-manifolds consists of $\SP^2\times \SP^2$ and the
blow-ups of $\C P^2$ at $k$, $0\leq k\leq 3$, points equipped with
the usual structures of a symplectic toric manifold (see
Theorem~\ref{thm-quantum-hom-sympl-toric-4-mfd-is-semi-simple-INTRO}
and Section~\ref{sec-appendix-quantum-hom-computations} below).
Another class of examples of symplectic manifolds admitting a
non-linear symplectic quasi-state was discovered recently by
D.~McDuff \cite{McD-private}: it includes a one point blow-up of
any non-uniruled closed symplectic manifold (for instance, of the
standard symplectic four-torus $\T^4$), see
Section~\ref{sec-fromto} below for more details. However we do not
know at the moment whether a non-linear symplectic quasi-state
exists on $\T^4$ (without blowing up)! In fact, there is only one
currently known way of constructing non-linear symplectic
quasi-states on symplectic manifolds of dimension higher than $2$.
It is based on the Floer theory for periodic orbits of Hamiltonian
flows and will be discussed in more detail below. Beyond it the
question of existence of a non-linear symplectic quasi-state on a
given symplectic manifold of dimension greater than two is
completely open.

Let us mention that on two-dimensional manifolds symplectic
quasi-states exist in abundance. They are supplied by the theory of
topological quasi-states, a beautiful recently emerged branch of
functional analysis whose foundations were laid by Aarnes in
\cite{Aar}.

Non-linear symplectic quasi-states are interesting for several
reasons.

The first reason is related to a mathematical model of quantum
mechanics which goes back to von Neumann \cite{von Neumann}. The
basic objects in this model are {\it observables} and {\it states}.
Observables are Hermitian operators on a complex Hilbert space $V$.
They form a real Lie algebra $\cA_q$ ($q$ for quantum) with the Lie
bracket given by $[A,B]_{\hbar} = \frac{i}{\hbar}(AB-BA)$, where
$\hbar$ is the Planck constant. States are real-valued linear
functionals on $\cA_q$ that are positive (i.e take non-negative
values on non-negative Hermitian operators) and normalized (that is
equal $1$ on the identity operator). Observables represent physical
quantities associated to a physical system such as energy, position,
momentum etc. while the value of a state on an observable represents
the mathematical expectation for the value of such a physical
quantity in the given state of the physical system.

A number of prominent physicists disagreed with the linearity
axiom put by von Neumann in the definition of a state. Their
reasoning was that the additivity of a state on two observables
$A$,$B$ makes sense {\it a priori} only if the corresponding
physical quantities are simultaneously measurable, meaning that
$A$ and $B$ commute: $[A,B]_{\hbar} = 0$. However in 1957 Gleason
\cite{Gleason} proved a remarkable theorem implying, under the
assumption ${\rm dim}\, V\geq 3$, that even if one considers a
positive functional on $\cA_q$ which is only required to be linear
on any pair of commuting observables then it has to be globally
linear. In other words, each "quantum quasi-state" has to be a
state.

In the analogous mathematical model of classical mechanics the
algebra $\cA_c$ of observables ($c$ for classical) is the space of
continuous functions $C(M)$ on a symplectic manifold $(M,\omega)$.
The Lie bracket is defined as the Poisson bracket on the dense
subspace $C^{\infty}(M) \subset C(M)$. States are linear
functionals on $\cA_c$ that are positive (i.e. take non-negative
values on non-negative functions) and normalized (equal $1$ on the
constant function $1$). Thus the existence of non-linear
symplectic quasi-states exhibits what can be called an
"anti-Gleason" phenomenon in classical mechanics: not every
"classical quasi-state" is a state \cite{EP-qst}, \cite{EPZ}.

The\, second\, point\, of\, interest\, in\, symplectic\,
quasi-states\, is\, related\, to\, the\, $C^0$-rigidity of the
Poisson brackets discovered by Cardin and Viterbo
\cite{Cardin-Viterbo}: while the definition of the Poisson
brackets $\{ F,G\}$ involves the first derivatives of $F$ and $G$,
it turns out that when the function $\{ F,G\}$ is not identically
zero its $C^0$-norm cannot be made arbitrarily small by
arbitrarily $C^0$-small perturbations of $F$ and $G$. Symplectic
quasi-states can be used to give effective quantitative
expressions for this and similar rigidity phenomena -- see
\cite{EPZ} for details.

Thirdly, applications of symplectic quasi-states include results
on rigidity of intersections -- a well-known phe\-no\-me\-non in
symplectic to\-po\-lo\-gy, meaning that a subset of a symplectic
manifold $(M,\omega)$ cannot be completely displaced from another
subset by a symplectic isotopy while it is possible to do it by a
smooth one. In fact, in all known applications of this sort it
suffices to have a functional on $C(M)$ with weaker properties,
called a {\it partial symplectic quasi-state}, whose definition
won't be given here -- see \cite{EP-qst}, \cite{EP-rigid} for
details. However in case the symplectic quasi-state $\zeta$ comes
from the Floer theory construction the reason it yields the
rigidity of intersections becomes really basic and simple
\cite{EP-qst}: by a generalized Riesz representation theorem of
Aarnes \cite{Aar} to each symplectic quasi-state $\zeta$ one can
associate a measure-type set-function $\tau_\zeta$, called a {\it
quasi-measure} (or a {\it topological measure}), which is finitely
additive (but not finitely subadditive unless $\zeta$ is linear).
Moreover, such a symplectic quasi-state $\zeta$ is invariant with
respect to the natural action of the group\footnote{The group
$\Symp_0 (M)$ is the identity component of the group $\Symp\, (M)$
of symplectomorphisms of $M$.} $\Symp_0 (M)$ on $C(M)$
("Invariance property") and therefore $\tau_\zeta$ is invariant
under symplectic isotopies. Hence if the sum of the quasi-measures
of two subsets of $M$ is greater than the quasi-measure of the
whole manifold $M$ none of these subsets can be displaced from the
other by a symplectic isotopy. In order to compute the
quasi-measures one may use in certain cases the following
important property of $\zeta$: it vanishes on functions with
sufficiently small support. More precisely, we call a set
$U\subset M$ {\it displaceable} if there exists a Hamiltonian
symplectomorphism\footnote{See e.g. \cite{McD-Sal-intro},
\cite{Pol-book} for preliminaries on Hamiltonian
symplectomorphisms.} $\psi$ such that $\psi (U) \cap
\text{Closure}\, {U} =\emptyset$. Then $\zeta (F) = 0$ for any
$F\in C(M)$ such that ${\rm supp}\, F$ is displaceable ("Vanishing
property").

\medskip
\noindent{\sc Calabi quasi-morphisms:}  Fourthly, any symplectic
quasi-state $\zeta$ coming from the Floer theory construction
mentioned above is an "infinitesimal version" of a function $\mu:
\tHam (M)\to\R$ on the universal cover $\tHam (M)$ of the group
$\Ham (M)$ of Hamiltonian symplectomorphisms of $M$ whose Lie
algebra can be identified with $C^\infty (M)/\R$. We shall need the
following notation: Given an (in general, time-dependent)
Hamiltonian $F=F_t$, $0\leq t\leq 1$, on $M$ denote by $\phi_F$ an
element of $\tHam (M)$ represented by the time-1 Hamiltonian flow
generated by $F$ and viewed as an identity-based path in $\Ham (M)$.
With this notation, for every $F \in C^{\infty}(M)$
\begin{equation}
\label{eqn-zeta-mu} \zeta (F) = \frac{\int_M F\omega^n}{{\rm vol}\,
(M)} - \frac{\mu (\phi_F)}{{\rm vol}\, (M)},
\end{equation}
where  ${\rm vol}\, (M) = \int_M \omega^n$ is the symplectic volume
of $M$.

 The function $\mu$ is of interest on its own merit -- it is
a {\it homogeneous quasi-morphism}, i.e. it satisfies the following
properties:

\medskip
\noindent \underline{\it Quasi-additivity:} There exists $C>0$,
which depends only on $\mu$, so that
\[
|\mu (\phi\psi) - \mu (\phi) - \mu(\psi)|\leq C \ \ {\rm for}\
{\it all}\ {\rm elements}\ \phi, \psi\in\tHam (M).
\]

\noindent \underline{\it Homogeneity:} $\mu (\phi^m) = m\mu
(\phi)$ for each $\phi\in\tHam (M)$ and each $m\in\Z$.

\medskip
Homogeneous quasi-morphisms on groups have become increasingly
important objects in various fields of mathematics (see e.g.
\cite{Ghys-ICM}),
especially in the case when a group does not
admit a non-trivial real-valued homomorphism and thus homogeneous
quasi-morphisms are, in a sense, the closest approximation to a
homomorphism one can hope to get on such a group. This is
precisely the case for $\tHam (M)$ which is a perfect group by
Banyaga's theorem \cite{Ban} and therefore does not admit a
non-trivial real-valued homomorphism. Moreover, the homogeneous
quasi-morphism $\mu$, constructed by means of the Floer theory,
has an additional remarkable property: it "patches up" a family of
homomorphisms on a certain class of subgroups of $\tHam (M)$ in
the following sense. For any open $U \subset M$, $U \neq M$,
define a subgroup $\tHam (U)$ of $\tHam (M)$ as
$$\tHam (U) :=
\{\, \phi_F\in \tHam (M)\ |\ \text{supp} (F_t) \subset U\; \text{for}\; \text{all}\; t\}.$$
Assuming that $F_t$ is supported in $U$ for every $t$, the formula
$$\phi_F \mapsto \int_0^1 dt \int_M F_t \omega^n$$
yields a correctly defined function $\Ca_U : \tHam (U) \to \R$
which is moreover a homomorphism, called the {\it Calabi
homomorphism} (see \cite{Ban,McD-Sal-intro}). Then $\mu$ satisfies
the following property corresponding on the infinitesimal level to
the vanishing property of $\zeta$:

\medskip
\noindent \underline{\it Calabi property:} $\left.
\mu\right|_{\widetilde{Ham}\, (U)} = \Ca_U$ for any displaceable open
$U\subset M$.

A homogeneous quasi-morphism on $\tHam (M)$ satisfying the Calabi
property is called a {\it Calabi quasi-morphism}. In particular,
we can now state

\begin{thm}
\label{thm-quantum-hom-sympl-toric-4-mfd-is-semi-simple-INTRO}
Assume $M$ is a symplectic toric Fano $4$-manifold with a rational
symplectic structure. Then $M$ admits a Calabi quasi-morphism $\mu$
and a symplectic quasi-state $\zeta$ which satisfies the vanishing
and the invariance properties, with $\mu$ and $\zeta$ being related
to each other by (\ref{eqn-zeta-mu}).
\end{thm}

For applications of the Calabi quasi-morphisms coming from the
Floer theory  to the study of the Hofer geometry on  $\Ham (M)$
and the algebraic structure of $\tHam (M)$ and $\Ham (M)$ see
\cite{EP-qmm}. Let us mention also that Py constructed Calabi
quasi-morphisms for closed surfaces of genus $\geq 1$ using
completely different ideas, see \cite{Py1,Py2}.

Let us mention that in \cite{EP-qst} we first introduced the Calabi
quasi-morphism $\mu$, then defined the symplectic quasi-state
$\zeta$ on the space $C^{\infty}(M)$ of smooth functions by formula
\eqref{eqn-zeta-mu}, and afterwards extended it to the space $C^(M)$
of continuous functions by the continuity in the uniform norm. We
give a concise definition of $\zeta$ in formula \eqref{eqn-def-zeta}
below.

\section{Quantum homology}\label{sec-prel}

Abusing the notation, we will write $\omega (A)$, $c_1 (A)$ for
the results of evaluation the cohomology classes $[\omega]$ and
$c_1 (M)$ on $A\in H_2 (M; \Z)$. Set
$$\bar{\pi}_2 (M) := H_2^S (M) / \sim,$$
where by definition
$$A \sim B\  {\rm iff}\  \omega (A) = \omega (B) \ {\rm and}\ c_1 (A) = c_1 (B).$$
 Denote by $\Gamma(M,\omega) := [\omega] (H_2^S
(M))\subset \R$ the subgroup of periods of the symplectic form on
$M$ on spherical homology classes.

\medskip
\noindent{\sc Novikov ring:} Let $\Gamma \subset \R$ be a
countable subgroup (with respect to the addition). Let $s,q$ be
formal variables. Let $\cF$ be a basic field of coefficients (in
applications, one usually takes $\cF=\C$ or $\cF=\Z_2$). Define a
field $\cK_{\Gamma}$ whose elements are generalized Laurent series
in $s$ of the following form:
$$\cK_{\Gamma} := \bigg\{\ \sum_{\theta \in \Gamma} z_\theta s^\theta, \ z_\theta \in
\cF,\ \sharp\big\{ \theta > c\ |\ z_\theta\neq 0\big\} <\infty,\
\forall c\in\R\ \bigg\}\;.$$ Define a ring $\Lambda_{\Gamma} : =
\cK_{\Gamma} [q, q^{-1}]$ as the ring of polynomials in $q,
q^{-1}$ with coefficients in $\cK_{\Gamma}$.  We turn
$\Lambda_{\Gamma}$ into a graded ring by setting the degree of $s$
to be zero and the degree of $q$ to be $2$. Note that the grading
on $\Lambda_{\Gamma}$ takes only even values. The ring
$\Lambda_{\Gamma}$ serves as an abstract model of the Novikov ring
associated to a symplectic manifold.

By definition, the Novikov ring of a symplectic manifold
$(M,\omega)$ is $\Lambda_{\Gamma(M,\omega)}$. In what follows,
when $(M,\omega)$ is fixed, we abbreviate and write  $\Gamma$,
$\cK$ and $\Lambda$ instead of $\Gamma(M,\omega)$,
$\cK_{\Gamma(M,\omega)}$ and $\Lambda_{\Gamma(M,\omega)}$
respectively.

\medskip
\noindent {\sc Quantum homology ring - definition:} The quantum
homology ring $QH_* (M)$ is defined as follows. First, it is a
graded module over $\Lambda$ given by
$$QH_* (M) := H_* (M; \cF)\otimes_{\cF}  \Lambda,$$
with the grading defined by the gradings on $H_* (M; \cF)$ and
$\Lambda$:
$$ {\rm deg}\, (a\otimes zs^\theta q^k) := {\rm deg}\, (a) + 2k.$$
Second, and most important, the $\Lambda$-module $QH_* (M)$ is
equipped with a {\it quantum product}: given $a\in H_k (M; \cF)$,
$b\in H_l (M; \cF)$, their quantum product is a class $a * b\in
QH_{k+l - 2n} (M)$\footnote{Recall that $2n = {\rm dim}\, M$.},
defined by
$$ a*b = \sum_{A\in \bar{\pi}_2 (M)} (a*b)_A \otimes s^{-\omega
(A)} q^{-c_1 (A)},$$ where $(a*b)_A \in H_{k+l - 2n + 2c_1 (A)} (M;
\cF)$ is defined by the requirement
$$ (a*b)_A \circ c = GW_A^{\cF} (a,b,c) \ \forall c\in H_* (M; \cF).$$
Here  $\circ$ stands for the intersection index and $GW_A^{\cF}
(a,b,c)\in \cF$ denotes the Gromov-Witten invariant which, roughly
speaking counts the number of pseudo-holomorphic spheres in $M$ in
the class $A$ that meet cycles representing $a,b,c\in H_* (M;\cF)$
(see \cite{Ru-Ti}, \cite{Ru-Ti-1}, \cite{MS2} for the precise
definition).

Extending this definition by $\Lambda$-linearity to the whole
$QH_* (M)$ one gets a correctly defined graded-commutative
associative product operation $*$ on $QH_\ast (M)$ which is a
deformation of the classical $\cap$-product in singular homology
\cite{Liu}, \cite{MS2}, \cite{Ru-Ti}, \cite{Ru-Ti-1}, \cite{Wi}.
The {\it quantum homology algebra} $QH_* (M)$ is a ring whose
unity  is the fundamental class $[M]$ and which is a module of
finite rank over $\Lambda$. If $a, b\in QH_\ast (M)$ have graded
degrees $deg\, (a)$, $deg\, (b)$ then \[deg\, (a\ast b) = deg\,
(a) + deg\, (b) - 2n.
\]

Denote by $H_{ev} (M; \cF)$ the even-degree part of the singular
homology and by $QH_{ev} (M) = H_{ev} (M;\cF)\otimes_{\cF}
\Lambda$ the even-degree part of $QH_* (M)$. Then $QH_{ev} (M)$ is
a commutative subring of $QH_* (M)$ which is a module of finite
rank over $\Lambda$. We will identify $\Lambda$ with a subring of
$QH_{ev} (M)$ by $\lambda \mapsto [M]\otimes \lambda$.

\medskip
\noindent{\sc Semi-simple algebras:} Re\-call that a
com\-mu\-ta\-tive (finite-di\-men\-sio\-nal) algebra $Q$ over a
field $\cB$ is called {\it semi-simple} if it splits into a direct
sum of fields as follows: $Q = Q_1 \oplus...\oplus Q_d\;$, where

\begin{itemize}
\item{}each $Q_i \subset Q$ is a finite-dimensional linear
subspace over $\cB$; \item{} each $Q_i$ is a field with respect to
the induced ring structure; \item{} the multiplication in $Q$
respects the splitting:
$$(a_1,...,a_d)\cdot(b_1,...,b_d) = (a_1 b_1,...,a_d b_d).$$
\end{itemize}
A classical theorem of Wed\-der\-burn (see e.g.
\cite{VanDerWaerden}) implies that  the semi-simplici\-ty is
equivalent to the absence of nilpotents in the algebra.

\medskip
\noindent  Let $Q$ be a commutative algebra over a field
$\cB$. Assume the base field $\cB$ is a subfield of a larger field
$\bar{\cB}$. Then $Q$ extends to a $\bar{\cB}$-algebra $\bar{Q}$:
$$\bar{Q} := Q\otimes_\cB \bar{\cB}.$$
One can identify $Q$ with a $\cB$-subalgebra of $\bar{Q}$
(whenever the latter is viewed as a ${\cB}$-algebra).

\begin{prop}
\label{prop-semisimple-algebras-basics}\

\medskip
\noindent (A) Assume the field $\cB$ is of characteristic zero. Then
the $\cB$-algebra $Q$ is semi-simple if and only if the
$\bar{\cB}$-algebra $\bar{Q}$ is semi-simple.

\medskip
\noindent (B) A tensor product of finite-dimensional semi-simple
algebras over a field of characteristic $0$ is again semi-simple.

\end{prop}

For the proof see e.g. \cite{VanDerWaerden}, \S 103. The "if" part
of (A) is trivial and does not require any assumptions on the
characteristic of $\cB$: if there were a nilpotent in $Q$ then it
would have been a nilpotent also in $\bar{Q}$.

A remark on part (A) of the proposition: Note that the fields in
the decomposition of the algebra, even their number, may change
and the unity of such a field may no longer be one after passing
from $Q$ to $\bar{Q}$.

Semi-simplicity is a partial case of the following more general
algebraic property which is relevant to our discussion -- see
Section~\ref{sec-fromto} below. Namely, we say a commutative
finite-dimensional algebra $Q$ over a field $\cB$ {\it contains a
field as a direct summand} if it decomposes as a $\cB$-algebra
into {\it a direct sum of $\cB$-algebras}
\begin{equation}
\label{eq-McDuff} Q =Q_1 \oplus Q_2\;,
\;\;\text{where}\;\; \;Q_1 \;\;\;\text{is a
field}\;.\end{equation} Let us emphasize that no assumption on
$Q_2$ is imposed.

The next proposition generalizes
Proposition~\ref{prop-semisimple-algebras-basics}. Though it is
not used in the sequel of the paper, it highlights an important
feature of the property of containing a field as a direct summand:
this property turns out to be robust with respect to extensions of
the base fields. This feature is highly desirable in the context
of quantum homology where one faces the necessity to choose the
base field depending upon the circumstances. As before, assume
$\cB$ is a subfield of a larger field $\bar{\cB}$ and extend $Q$
to a $\bar{\cB}$-algebra $\bar{Q}$: $\bar{Q} := Q\otimes_\cB
\bar{\cB}$.

\begin{prop}
\label{prop-direct-summand}\

\medskip
\noindent (A) Assume the field $\cB$ is of characteristic zero. Then
the $\cB$-algebra $Q$ contains a field as a direct summand if and only if the
$\bar{\cB}$-algebra $\bar{Q}$ does.

\medskip
\noindent (B) Assume two finite-dimensional commutative algebras
over a field of characteristic $0$ both contain fields as direct
summands. Then their tensor product also has this property.

\end{prop}

\medskip
\noindent The "only if" part of (A) and part (B) follow readily
from Proposition~\ref{prop-semisimple-algebras-basics}. The "if"
part of (A) is  less trivial and will be proved now for the sake
of completeness.

\bigskip
{\sc Proof.}

\medskip
\noindent{\sc Preliminaries:} First of all, it suffices to prove
the theorem for algebraically closed $\bcB$: indeed, by part (A)
of Proposition 2.1, if $Q \otimes_\cB \bcB$ contains a field as a
direct summand, the same holds for $Q \otimes_\cB \bcB'$ where
$\bcB'$ is the algebraic closure of $\bcB$. In what follows we
assume that $\bcB$ is algebraically closed.

We shall use notation $A \oplus B$ for the direct sum of {\it
algebras} (over the same base field) and $A \dot{+} B$ for the
direct sum of {\it modules} (also over the same field).

Given extensions $\cB \subset E \subset \bcB$ of fields, we shall
say that an element $t \in Q \otimes_\cB \bcB$ is {\it defined over}
$E$ if it can be written as $t = \sum q_j \otimes \theta_j$ with
$q_j \in Q, \theta_j \in E$.

For $q \in Q$ we write $\widehat{q} := q \otimes 1_{\bcB}$, where
$1_{\bcB}$ is the unity of $\bcB$.

Every semi-simple commutative finite-dimensional algebra over an
algebraically closed field admits a unique  decomposition into
one-di\-men\-sio\-nal fields. The correspondent unities (called
{\it basic idempotents}) sum up to the unity of the algebra (we
refer to this as to {\it the partition of unity}). By the
Krull-Schmidt theorem, the basic idempotents are unique up to
permutation.

\medskip
\noindent {\sc Starting the proof:} By the Wedderburn principal
theorem (see e.g. \cite{Faith}, Thm 13.18)
$$Q = (F_1 \oplus ... \oplus F_n) \dot{+} Y\;,$$
where $Y$ is the radical \footnote{For finite-dimensional
commutative algebras the radical is the set of all nilpotent
elements.} and all $F_i$ are fields. Consider any monomorphism of
$F_i$ into $\bcB$ and denote by $E_i \subset \bcB$ the normal
closure of the image. The field $E_i$ does not depend on the choice
of a monomorphism. A concrete realization of $E_i$ is as follows.
Take any primitive element $a_i \in F_i$ and denote by $P_i$ its
minimal polynomial. Consider the factorization of $P_i$ in $\bcB$
into linear factors:
$$P_i(u) = (u-\mu_{i1})\cdot...\cdot(u-\mu_{ir_i})\;.$$
All $\mu_j$'s are pair-wise distinct (since $F_i$ is a field). Then
$$E_i = \cB (\mu_{i1},...,\mu_{ir_i})\;.$$

\medskip
\noindent
\begin{lemma}\label{lem-1} The basic idempotents
of $F_i \otimes_\cB \bcB$ are defined over $E_i$ for all $i=1, ... ,
n$.
\end{lemma}

{\sc Proof of the lemma.} We use the notation introduced above,
and omit the lower index $i$. Given $j \in \{1;...;r\}$, introduce
the polynomial $t_j(u)$ of degree $r -1$ which vanishes at all
$\mu_m$ for $m \neq j$ and equals $1$ at $\mu_j$. One immediately
checks that its coefficients lie in $E$. Note that the basic
idempotents of $F \otimes_\cB \bcB$ are given by $t_j(a)$: In
order to see this use that $F = \cB [a] = \cB [u]/(P)$ and hence
$F \otimes_\cB \bcB = \bcB [u]/(P)$. The lemma follows. \Qed

\medskip

Denote by $G_i$ the Galois group of $E_i$ over $\cB$. We extend its
action to $Q\otimes_{\cB} E_i$ by the change of coefficients:
$$g(q \otimes \theta) = q \otimes g(\theta) \;\; \forall q \in Q,\;
\theta \in E_i,\;g \in G_i\;.$$

\begin{lemma}\label{lem-2} The Galois group $G_i$ acts transitively
on  basic idempotents of $F_i \otimes_\cB \bcB$ for all $i=1, ... ,
n$.
\end{lemma}

{\sc Proof of the lemma.} Again, we omit the lower index $i$. Take
any basic idempotent, say $e_1$ and denote by $\mathcal O$ its
orbit under the action of $G$. Then
\begin{equation}\label{eq-e}
h:= \sum_{f \in \mathcal{O}} f
\end{equation}
is an idempotent {\it invariant under the action of} $G$. Indeed,
the invariance is obvious and in order to prove that $h$ is an
idempotent, note that, according to Lemma~\ref{lem-1}, each $g\in G$
preserves the set of basic idempotents of $F \otimes_{\cB} \bcB$ and
therefore $h$ is the sum of the basic idempotents lying in the orbit
of the action of $G$ on that set, which implies that $h$ itself is
an idempotent.

Moreover, $h$ is defined over $\cB$. Therefore it lies in $F
\otimes_\cB \cB$. The latter is isomorphic to $F$, which was assumed
to be a field. Thus $h = \widehat{1_F}$. We conclude that formula
\eqref{eq-e} gives the partition of unity in $F \otimes_{\cB}
 \bcB$. Thus the sum at the right hand side of the
formula includes {\it all} the basic idempotents, and the lemma
follows. \Qed

\medskip
\noindent {\sc Exploring the direct summand:} Let $H\subset
Q\otimes_\cB \bcB$ be the field (the direct summand) from the
assumption of the theorem. Denote by $e$ its unity.

\medskip
\noindent
\begin{lemma}\label{lem-vsp-idemp}
Element $e$ coincides with one of the basic idempotents of some $F_i
\otimes_\cB \bcB$.
\end{lemma}

{\sc Proof of the lemma.} Since
$$Q\otimes_\cB \bcB = (F_1 \otimes_\cB \bcB \oplus ... \oplus F_n \otimes_\cB
\bcB) \dot{+} Y \otimes_\cB \bcB$$ we can decompose
$$e = (x_1 \oplus ... \oplus x_n) \dot{+} y$$
with $x_i \in F_i \otimes_\cB \bcB, y \in Y \otimes_\cB \bcB$.

Since $e^2=e$ we have
$$(x_1^2 \oplus ... \oplus x_n^2) \dot{+} (2x_1y + ... + 2x_ny +
y^2) = (x_1 \oplus ... \oplus x_n) \dot{+} y\;.$$ Thus
\begin{equation}\label{eq-vsp-x}
x_i^2 =x_i \;\; \forall i =1,...,n
\end{equation}
and
\begin{equation}\label{eq-vsp-y}
2x_1y + ... + 2x_ny + y^2=y\;.
\end{equation}
Multiplying \eqref{eq-vsp-y} by $x_i$ we have $x_iy + x_iy^2 = 0$.
Rewrite this as $x_iy + (x_iy)^2 =0$. Since $x_iy$ is a nilpotent,
this yields $x_iy=0$. Substituting this into \eqref{eq-vsp-y} we get
$y^2 = y$, and since $y$ is nilpotent we have $y=0$.

Let us write $x_i$ as a linear combination of the basic
idempotents $e_{ij}$, $j=1,...,r_i$, of $F_i \otimes_\cB \bcB$ and
substitute into \eqref{eq-vsp-x}. We get that
\begin{equation}\label{eq-vsp-sum}
e = \sum_{i=1}^n \sum_{j \in S_i} e_{ij}
\end{equation}
for some subsets $S_i \subset \{1;...;r_i \}$. We claim that the
sum in the right hand side of \eqref{eq-vsp-sum} consists of one
term only: Indeed, assume on the contrary that there are at least
two distinct terms, say $e_{ab}$ and $e_{cd}$. Then $e_{ab} =
e_{ab}e$ and $e_{cd}= e_{cd}e$ and therefore both $e_{ab}$ and
$e_{cd}$ lie in $H$ (here we use that $H$ is a direct summand in
the sense of algebras). But $e_{ab}\cdot e_{cd}=0$, and we get a
contradiction with the fact that the field $H$ has no divisors of
zero. The claim, and thus the lemma, follow. \Qed

\medskip
\noindent{\sc The end of the proof of the proposition:} In view of
Lem\-ma~\ref{lem-vsp-idemp}, we can assume without loss of
generality that the unity element $e$ of $H$ equals $e_{11}$.

Since $H$ is a direct summand and a field, $\widehat{y} \cdot e_{11}
= 0$ for every $y$ in the radical $Y$. We claim that $\widehat{y}
\cdot e_{1j} = 0$ for all $j=1,...,r_1$. Indeed, for given $j$ by
Lemma \ref{lem-2} there is an element, say, $g$ of the Galois group
$G_1$ with $g(e_1)=e_j$. Note that $G_1$ fixes $\widehat{q}$ for
every $q \in Q$. Thus
$$\widehat{y} \cdot e_{1j} = g(\widehat{y} \cdot
e_{11}) = 0\;,$$ and the claim follows. Since
$$e_{11} + ... + e_{1r_1} = \widehat{1_{F_1}}$$
we have that $$\widehat{y\cdot 1_{F_1}} =\widehat{y} \cdot
\widehat{1_{F_1}}=0 \;\; \forall y \in Y\;.$$ Thus $Y\cdot F_1 =0$
in $Q$, which yields the decomposition
$$Q= F_1 \oplus (F_2 \oplus ... \oplus F_n \dot{+} Y)\;.$$
We conclude that $F_1$ is the desired field summand of $Q$. This
completes the proof. \Qed

\section{From quantum homology to quasi-states}\label{sec-fromto}

\noindent{\sc Spectral invariants:}  Let us now discuss in more
detail the construction of $\zeta$ via Hamiltonian Floer theory.
The construction involves {\it spectral numbers} $c (a, H)$ and
$c(a, \phi)$ associated by means of certain minimax-type procedure
in the Floer theory to any non-zero quantum homology class $a$,
any (time-dependent) Hamiltonian $H$ on $M$  \cite{Oh-spectral}
and any element $\phi \in \tHam(M,\omega)$ (see also
\cite{Viterbo, Oh1, Schwarz, Oh2} for earlier versions of this
theory). We refer to \cite{EP-rigid} for a detailed review.

In what follows we assume that all the symplectic manifolds we
deal with belong to the class $\cS$ of closed symplectic manifolds
for which the spectral invariants are well-defined and enjoy the
standard list of properties (see e.g. \cite[Theorem 12.4.4]{MS2}).
For instance, $\cS$ contains all spherically monotone manifolds,
that is symplectic manifolds $(M,\omega)$ such that $\left.
[\omega]\right|_{H_2^S (M)} =\kappa \left. c_1 \right|_{H_2^S
(M)}$ for some $\kappa
>0$, as well as manifolds with $\left. [\omega]\right|_{H_2^S (M)} = 0$.
Furthermore, $\cS$ contains all symplectic manifolds $M^{2n}$ for
which, on one hand, either $c_1=0$ or the minimal Chern number
(i.e. the positive generator of $c_1 (H_2^S (M)) \subset \Z$) is
at least $n-1$ and, on the other hand, $[\omega]$ is rational. This
includes all symplectic $4$-manifolds with a rational symplectic
structure. The general belief is that the class $\cS$ includes
{\bf all} symplectic manifolds\footnote{According to the recent
preprint of M.Usher \cite{Usher} and the previous result of
Y.-G.Oh \cite{Oh-nondeg-spectrality} on the "non-degenerate
spectrality" property of spectral numbers, the rationality
assumption on $[\omega]$ is not needed for the further results in
this paper.}.

Denote by $QH_{2n}(M)$ the graded component of degree $2n$ of the
quantum homology algebra  $QH_* (M)$. It is an algebra over the
field $\cK$. We will say that $QH_{2n} (M)$ is  semi-simple, if it
is semi-simple as a $\cK$-algebra. Elementary grading
considerations show that any non-zero idempotent in $QH_* (M)$ has
to lie in $QH_{2n} (M)$. Given an non-zero idempotent $a\in
QH_{2n} (M)$, a time-independent Hamiltonian $H: M\to \R$ and an
element $\phi\in \tHam (M)$, define
\begin{equation}
\label{eqn-def-zeta} \zeta (H): = \lim_{l \to +\infty}\; \frac{c
(a,lH)}{l},
\end{equation}
and
\begin{equation}
\label{eqn-def-mu} \mu (\phi): = - {\hbox{\rm vol}\,
(M)}\cdot\lim_{l \to +\infty}\; \frac{c (a,\phi^l)}{l},
\end{equation}
where ${\hbox{\rm vol}\, (M)} := \int_M \omega^n$. For manifolds
from class $\cS$  formula \eqref{eqn-def-zeta} yields a functional
on $C^\infty (M)$ which can be extended to $C(M)$. The link
between existence of symplectic quasi-states and Calabi
quasi-morphisms on the one hand, and semi-simplicity of quantum
homology on the other hand, is given by the next theorem.

\begin{thm}
\label{thm-qstate-non-monotone} Assume that $QH_{2n} (M)$ is
semi-simple and decomposes into a direct sum of fields: $QH_{2n}
(M) = Q_1\oplus\ldots\oplus Q_d$. Assume the idempotent $a$
appearing in the definitions of $\zeta$ and $\mu$ (see equations
\eqref{eqn-def-zeta}, \eqref{eqn-def-mu}) is the unit element in
some $Q_i$. Then $\mu$ is a Calabi quasi-morphism and $\zeta$ is a
symplectic quasi-state satisfying the vanishing and the invariance
properties.
\end{thm}

\medskip \noindent{\sc The role of semi-simplicity:} The example
$M= \T^2$ \cite{EP-qst} shows that the semi-simplicity assumption
cannot be completely omitted -- without it $\zeta$ may turn out to
be only a partial symplectic quasi-state (a weaker object already
mentioned above) and not a genuine one. On the other hand, as it
was pointed out to us by D.~McDuff \cite{McD-private}, it suffices
to assume that the $\cK$-algebra $QH_{2n} (M)$ contains a field as
a direct summand (see \eqref{eq-McDuff} above for the definition).
Indeed, take the idempotent $a$ appearing in the definitions of
$\zeta$ and $\mu$ (see equations \eqref{eqn-def-zeta},
\eqref{eqn-def-mu}) to be the unit element in the field appearing
as a direct summand of $QH_{2n} (M)$. Then $\mu$ is a Calabi
quasi-morphism and $\zeta$ is a symplectic quasi-state satisfying
the vanishing and the invariance properties. The proof of this
statement repeats {\it verbatim} the one of
Theorem~\ref{thm-qstate-non-monotone}. The results of McDuff
\cite{McD-Uniruled} yield a remarkable class of symplectic
manifolds whose quantum homology, while in general being
non-semi-simple, admits a decomposition \eqref{eq-McDuff}: These
are one point blow-ups of {\it non-uniruled} symplectic manifolds
(for instance, of the standard symplectic torus $\T^4$ or of any
other symplectic manifold of dimension greater than 2 with $\left.
[\omega]\right|_{\pi_2 (M)} = 0$) -- see
Section~\ref{sec-mcduff-examples} for precise definitions and
statements.

\medskip \noindent{\sc Various versions of semi-simplicity:}
Originally the construction of $\mu$ and $\zeta$ was carried out
in \cite{EP-qmm,EP-qst}  for a spherically monotone $M$ under the
assumption that the even-degree quantum homology of $M$ is
semi-simple as a commutative algebra over a field of Laurent-type
series in one variable with complex coefficients (this field can
be identified with the Novikov ring in the monotone case).
Theorem~\ref{thm-semisimple-monotone-case} below states that this
assumption is equivalent to semi-simplicity of $QH_{2n}(M)$ in the
sense of the present paper.

An immediate difficulty with extending the result of \cite{EP-qmm}
to non-monotone manifolds lies in the fact that in the
non-monotone case the No\-vi\-kov ring is no longer a field. In
\cite{Ostr-qmm} Ostrover bypasses this difficulty by working only
with the graded part of degree $2n={\rm dim}\, M$ of the quantum
homology which is already an algebra over a certain field
contained in the Novikov ring -- in this way he proves the
existence of a Calabi quasi-morphism and a non-linear symplectic
quasi-state on $\SP^2\times \SP^2$ and $\C P^1\sharp \overline{\C
P^1}$ equipped with non-monotone symplectic forms. In addition to
the semi-simplicity Ostrover assumes that the minimal Chern number
$N$ of $M$ divides $n$. The definition of quantum homology
introduced in Section~\ref{sec-prel} above is an algebraic
extension of the one used by Ostrover: our formal variable $q$,
which is responsible for the Chern class, is the root of degree
$N$ of the corresponding variable in Ostrover's setting. This
minor modification is nevertheless useful: Ostrover's argument
(see \cite[Section 5]{Ostr-qmm})  can be applied {\it verbatim} to
the proof of Theorem \ref{thm-qstate-non-monotone} without the
assumption $N\; |\; n$. Note also that, when we work with
$\C$-coefficients, algebraic extensions respect semi-simplicity
(see Proposition \ref{prop-semisimple-algebras-basics}(A) above),
and hence semi-simplicity in our setting is equivalent to
semi-simplicity in Ostrover's setting.

\begin{rem}
{\rm There is another notion of semi-simple quantum (co)homology
that has been extensively studied by various authors -- see e.g.
\cite{Bayer}, \cite{Bayer-Manin}, \cite{Dubrovin}, \cite{Manin},
\cite{Tian-Xu} and the references therein. We do not know whether
in any of the examples the semi-simplicity in our sense can be
formally deduced from the semi-simplicity in the other sense. }
\end{rem}

\fbox{\parbox{0.8 \linewidth}{ \centerline{{\bf From now on we
work with the basic field $\cF=\C$.}} }}

\bigskip
\bigskip
\noindent {\sc Application to symplectic toric Fano 4-manifolds:}

\begin{thm}
\label{thm-quantum-hom-sympl-toric-4-mfd-is-semi-simple}  If $M$
is a symplectic toric Fano $4$-manifold with a rational symplectic
form, the algebra $QH_4 (M)$ is semi-simple.
\end{thm}

 {\sc Question.} Is it true that $QH_{2n}(M)$ is
semi-simple for any symplectic toric Fano manifold $M^{2n}$?

\medskip
\noindent The proof of
Theorem~\ref{thm-quantum-hom-sympl-toric-4-mfd-is-semi-simple}
proceeds as follows.  We inspect the list of all symplectic toric
Fano $4$-manifolds -- up to rescaling the symplectic form by a
constant factor, there are only five of them: $\SP^2\times \SP^2$,
$\C P^2$ and the blow-ups of $\C P^2$ at $1,2,3$ points. For $\C
P^2$ semi-simplicity is elementary (see e.g. \cite{EP-qmm}), and
for $\SP^2\times \SP^2$ and $\C P^2$ blown up at one point it was
established by Y.Ostrover in \cite{Ostr-qmm}. The cases of the
blow-ups of $\C P^2$ at $2$ or $3$ points are new. For the
calculations of the quantum homology, we use Batyrev's algorithm
(see \cite{Bat}, cf. \cite{MS2}) which involves the
combinatorics/geometry of the moment polytope, see
Section~\ref{sec-appendix-quantum-hom-computations} for the
details.

\medskip
\noindent
Theorem~\ref{thm-quantum-hom-sympl-toric-4-mfd-is-semi-simple-INTRO}
follows immediately from Theorem~\ref{thm-qstate-non-monotone} and
Theorem~\ref{thm-quantum-hom-sympl-toric-4-mfd-is-semi-simple}.

\medskip
\noindent In addition to symplectic toric Fano $4$-manifolds with
a rational symplectic structure other symplectic manifolds with
semi-simple quantum homology include complex projective spaces and
complex Grassmannians \cite{EP-qmm}. We show (see
Theorem~\ref{thm-product} below) that by taking products one can
construct more examples of manifolds with semi-simple quantum
homology. Moreover, by taking direct products of manifolds whose
quantum homology contains a field as a direct summand one can
construct more manifolds with the same property (see
Theorem~\ref{thm-field-direct-summand-product} below).

\section{Semi-simplicity of $QH_4 (M)$ for symplectic toric Fano $M^4$}
\label{sec-appendix-quantum-hom-computations}

\subsection{Classification of symplectic toric Fano $4$-manifolds}
\label{subsec-monotone-4-dim-toric-mfds-quantum-hom-computations}

\ \ Any closed symplectic toric manifold $(M,\omega)$ can be
equipped canonically with a complex structure $J$ invariant under
the torus action and also compatible with $\omega$ (meaning that
$\langle\cdot,\cdot\rangle := \omega (\cdot, J\cdot)$ is a
$J$-invariant Riemannian metric). Such a $J$ is unique up to an
equivariant biholomorphism -- see e.g. \cite{LT}, Section 9. We
say that the symplectic toric manifold $(M,\omega)$ is {\it Fano}
if $(M,J)$ is Fano as a complex variety. The symplectic toric Fano
$4$-manifolds can be completely described as follows:

\begin{prop}
\label{prop-Fano-4-mfds-classification} If $M$ is a Fano
symplectic toric $4$-manifold then $M$ is equivariantly
symplectomorphic to one of the following toric manifolds:\break
$\SP^2\times \SP^2$ or  the blow-up of $\C P^2$ at $k$, $0\leq
k\leq 3$, points, where the symplectic toric structures on
$\SP^2$, $\C P^2$ are the standard ones and the blow-ups are
equivariant (but the sizes of the blow-ups may vary).

\end{prop}

\begin{proof} Denote by $J$ the $\T^2$-invariant complex structure
on $M$ compatible with $\omega$. It is known (see e.g.
\cite[Appendix 1]{Guil}) that $(M,J)$ is a toric variety whose
fan, say $\Sigma$, is determined by the moment polytope $\Delta$
of the Hamiltonian $\T^2$-action as follows: To $i$-th edge of
$\Delta$ assign a perpendicular primitive integral vector $e_i$
pointing inside the polytope. View each $e_i$ as a vector in the
dual space $\R^{2*}$. The fan $\Sigma \subset (\R^{2})^*$ is
defined by the rays $\ell_i$ generated by $e_i$'s.

Since $(M,J)$ is Fano, the fan $\Sigma$ is spanned by the edges of
{\it a Fano polytope}, say $P$. Recall that a convex polytope in
$(\R^{2})^*$ is called Fano if it contains the origin in its
interior, its vertices lie in the integral lattice $\Z^2 \subset
(\R^{2})^*$ and each pair of consecutive vertices forms a basis of
$\Z^2$ (see \cite[p.304]{Ewa}). Since the vertices of $P$ are
primitive integer vectors lying in $\ell_i$'s, they coincide with
$e_i$'s. We conclude that the vectors $e_i$ form the set of
vertices of a Fano polytope. All such polytopes are classified --
there are precisely $5$ of them, up to the action of $GL (2,\Z)$
(see \cite{Klyachko-Voskr}, cf. \cite[p.192]{Ewa}), and their fans
define the complex toric surfaces listed in the proposition.
Looking at polygons $\Delta$ having $e_i$'s as inward normals, we
get that they correspond to the moment polygons of the standard
symplectic toric structures on these surfaces. Thus, by the
Delzant theorem \cite{Delz}, $(M,\omega)$ is equivariantly
symplectomorphic to one of these models. \end{proof}

Let us now prove
Theorem~\ref{thm-quantum-hom-sympl-toric-4-mfd-is-semi-simple}. The
semi-simplicity of the quantum homology algebra for $\SP^2\times
\SP^2$, $\C P^2$, $\C P^2\sharp \overline{\C P^2}$ was already
proved in \cite{EP-qmm}, \cite{Ostr-qmm}. Thus it remains to show
the claim only for the blow-ups of $\C P^2$ at $2$ and $3$ points.

\subsection{A recipe for computing
$QH_* (M)$ for a symplectic toric Fano $M^4$}

Note that a toric manifold has only even-dimensional integral
homology and therefore $QH_* (M) = QH_{ev} (M)$. The recipe for
computing the quantum homology of a symplectic toric Fano
$4$-manifold can be extracted from \cite{Bat}, \cite{MS2}. We shall
present it right now in a somewhat simplified form assuming that $M
\neq \C P^2$. Consider a moment polytope $\Delta$ of $M$. Let $l$
denote the number of its facets. To each facet $\Delta_i$,
$i=1,\ldots, l$, of $\Delta$ assign a perpendicular primitive
integral vector $e_i = (\alpha_i, \beta_i)$ pointing inside the
polytope and a two-dimensional homology class $u_i$ represented by
the pre-image of the facet under the moment map. Write the equation
for the line containing $\Delta_i$ as $\alpha_i x + \beta_i y =
\eta_i$, where $x,y$ are the coordinates on the space $\R^2$ where
$\Delta$ lies. The classes $u_i$ are additive generators of $H_2 (M;
\Z)$ and multiplicative generators of $QH_* (M)$. There are two
additive relations between these generators:
$$ \sum_{i=1}^l \alpha_i u_i = 0,\ \ \sum_{i=1}^l \beta_i u_i = 0.$$
The multiplicative relations can be described in the following
way. First of all, these relations correspond precisely to
primitive subsets $\{ i,j\}$ of $\{ 1,\ldots, l\}$ -- the
definition of a primitive subset in the case ${\rm dim}\, M=4$,
$M\neq \C P^2$, can be given as follows: a subset $I = \{ i,j\}$
of $\{ 1,\ldots, l\}$ is called {\it primitive} if the facets
$\Delta_i$ and $\Delta_j$ do not intersect. Given a primitive set
$I$, set $w:=\sum_{i\in I} e_i$. Let $\Sigma\subset (\R^2)^*$ be
the fan spanned by $e_i$'s (see the proof of
Proposition~\ref{prop-Fano-4-mfds-classification} above). Consider
the minimal cone in $\Sigma$ that contains $w$ and assume that the
cone is spanned by vectors $e_j$, $j\in J$ (if $w=0$ and the cone
is zero-dimensional we set $J=\emptyset$). One can show that the
set $J$ is always disjoint from $I$. If $w\neq 0$ then it can be
uniquely represented as $w=\sum_{j\in J} c_j e_j$, $c_j\in \N$.
Following \cite{Bat}, define an integral vector $d = (d_1, \ldots,
d_l)$ as follows: $d_k =1$ if $k\in I$, $d_k = -c_k$ if $k\in J$
and $d_k=0$ otherwise\footnote{Note that the definition of $d$ in
\cite{MS2}, p.414, does not necessarily lead to a unique $d$. In
order to get the uniqueness one should add that $d_\nu=0$ unless
$\nu\in J$ (with $J$ as defined here). Alternatively, one can drop
the uniqueness condition altogether and still get valid
multiplicative relations, though a bigger number of them. We thank
D.McDuff for clarifying this issue to us.}.
 Then the multiplicative relation in $QH_* (M)$
corresponding to $I$ is
$$ u_i * u_j =
s^{\sum_{k=1}^l d_k \eta_k  }q^{-\sum_{k=1}^l d_k} \prod_{m\notin
I} u_m^{-d_m}.$$

Finally we note that $QH_4 (M)$ is generated, as a subring of
$QH_* (M)$, by the elements $q u_i$.

In what follows we implement the recipe above for equivariant
blow-ups of $\C P^2$ at two and three points. We assume that the
standard symplectic form on $\C P^2$ is normalized in such a way
that the integral over the line equals $1$, and consider the
standard Hamiltonian torus action on $\C P^2$ whose moment polytope
is the triangle $\Pi$ with the vertices $(0,0),(1,0)$ and $(0,1)$ in
$\R^2$. These vertices correspond to the fixed points of the action
denoted respectively by $P_{00},P_{10}$ and $P_{01}$.

\subsection{The case of a blow-up of $\C P^2$ at 2 points}
\label{subsec-CP2-blown-up-at-2pts-quantum-hom-computations}

Let $M$ be obtained as a result of the equivariant blow-up of $\C
P^2$ at the fixed points $P_{10}$ and $P_{01}$ of the standard
action, so that the symplectic areas of the exceptional divisors
are equal to $1-\delta$ and $1-\epsilon$, where $\epsilon,\delta
\in (0;1)$ are rational numbers satisfying $\epsilon + \delta >1$.
The moment polytope $\Delta$ of the blown up torus action on $M$
is obtained from the triangle $\Pi$ by chopping off the corners
(see e.g. \cite[Section 3]{KK}) at the vertices $(1,0)$ and
$(0,1)$. Thus $\Delta$ is a pentagon with the vertices $(0, 0)$,
$(0, \epsilon)$, $(1-\epsilon, \epsilon)$, $(\delta, 1-\delta)$,
$(\delta, 0)$.

\begin{prop}
The $\cK$-algebra $QH_4 (M)$ is semi-simple.
\end{prop}

\begin{proof}

We number the facets of $\Delta$
as follows: the facet connecting $(0, \epsilon)$ and
$(1-\epsilon, \epsilon)$ is numbered by $1$, then we number the rest of the
faces moving clockwise from the first one. The vectors $e_i$
are the following: $e_1 = (0,-1)$, $e_2 = (-1,-1)$, $e_3=(-1,0)$,
$e_4 = (0,1)$, $e_5 = (1,0)$. The constants $\eta_i$ are as follows:
$\eta_1 = -\epsilon$, $\eta_2 = -1$, $\eta_3 = -\delta$,
$\eta_4 = \eta_5 = 0$.

The additive relations are:
$$ -u_2 - u_3 + u_5 = 0,\ \  -u_1 - u_2 + u_4 = 0.$$
We choose the basis of $H_2 (M;\Z)$ as $u_1, u_2, u_3$, and thus
$u_5 = u_2 + u_3$, \break $u_4 = u_1 + u_2$.

There are five primitive sets: $\{ 1, 3\}$, $\{ 1, 4\}$,
$\{ 2, 4\}$, $\{ 2, 5\}$, $\{ 3, 5\}$. Taking into account
the additive relations the corresponding multiplicative
relations for $QH_* (M)$ can be written as follows:
$$u_1 * u_3 = s^{-(\epsilon +\delta -1)} q^{-1}u_2,$$
$$u_1 * (u_1 + u_2) = s^{-\epsilon} q^{-2} [M],$$
$$u_2 * (u_1 + u_2) = s^{\delta-1} q^{-1} u_3,$$
$$u_2 * (u_2 + u_3) = s^{\epsilon-1} q^{-1} u_1,$$
$$u_3 * (u_2 + u_3) = s^{-\delta} q^{-2} [M].$$

To simplify the notation we shall drop the $*$ sign for the quantum
product, and write $1$ instead of $[M]$.

We express $u_2$ from the first relation via $u_1, u_3$, then
define $X:=qu_1$, $Y:= qu_3$ and substitute the result in the
remaining $4$ relations. It turns out that $X$, $Y$ are invertible
in $QH_* (M)$ and only $2$ of the $4$ relations are independent:
$$s^\epsilon X^2(1+s^{\epsilon+\delta-1}Y)=1,$$
$$s^\delta Y^2(1+s^{\epsilon+\delta-1}X)=1.$$

Factoring $\cK [X,Y]$ by these $2$ relations yields a complete
description of the $\cK$-algebra $QH_4 (M)$. Expressing $Y$ via
$X$ from the first relation and substituting the result in the
second shows yields the relation:
\begin{equation}\label{eq-rel}
s^{1-\epsilon}X^5 + (s^{2-2\epsilon-\delta}-1)X^4
-2s^{1-2\epsilon}X^3 - 2s^{2-3\epsilon-\delta}X^2 +
s^{1-3\epsilon}X + s^{2-4\epsilon-\delta} = 0.
\end{equation}
Denote the polynomial of $X$ in the left-hand side by
$Q_{\epsilon,\delta}(X;s)\in \cK[X]$. We have obtained that $QH_4
(M)\cong \cK [X,X^{-1}]/(Q_{\epsilon,\delta})$. Note that
$$Q_{\epsilon,\delta}(X;1)=f(X)\;,$$
where $f(X):= X^5 - 2X^3 - 2X^2 +X +1$. Denote by
$R_{\epsilon,\delta}(s)$ the resultant of
$Q_{\epsilon,\delta}(X;s)$ and its derivative with respect to
$X$-variable. Note that \break $R_{\epsilon,\delta}(1)=r$, where
$r$ is the resultant of $f$ and $f'$. A direct check shows that
$r\neq 0$, and hence $R_{\epsilon,\delta}(s) \neq 0\in \cK$. This
shows that $Q_{\epsilon,\delta}$ does not have a common root with
its derivative for any choice of $\epsilon$ and $\delta$, and thus
has no multiple roots in any extension of $\cK$. Thus the
$\cK$-algebra $\cK [X]/(Q_{\epsilon,\delta})$ is semi-simple. In
view of relation \eqref{eq-rel} $X^{-1}$ is a polynomial of $X$.
Thus ${\cK} [X,X^{-1}]/(Q_{\epsilon,\delta})=\cK
[X]/(Q_{\epsilon,\delta})$, and hence $QH_4(M)$ is semisimple.
This finishes the proof of the proposition. \end{proof}

\subsection{The case of a monotone blow-up of $\C P^2$ at 3 points}
\label{subsec-CP2-blown-up-at-3pts-quantum-hom-computations}

Let $M$ be obtained as a result of the equivariant blow-up of $\C
P^2$ at the fixed points $P_{00},P_{10}$ and $P_{01}$ of the
standard action, so that the symplectic areas of the exceptional
divisors are equal to $\alpha$, $1-\beta$ and $1-\gamma$, where
$\alpha,\beta,\gamma \in (0;1)$ are rational numbers satisfying
$\alpha < \gamma$, $\alpha < \beta$ and $\beta+\gamma >1$. The
moment polytope $\Delta$ of the blown up torus action on $M$ is
obtained from the triangle $\Pi$ by chopping off the corners (see
e.g. \cite[Section 3]{KK}) at its vertices. Thus $\Delta$ is a
hexagon with the vertices $(\alpha, 0)$, $(0,\alpha)$, $(0,
\gamma)$, $(1-\gamma,\gamma)$, $(\beta, 1-\beta)$, $(\beta, 0)$.

\begin{prop}
$QH_4 (M)$ is a semi-simple algebra over $\cK$.
\end{prop}

\begin{proof}

By technical reasons, it would be convenient to extend the
coefficients of the quantum homology from the field
$\cK:=\cK_{\Gamma}$, where $\Gamma$ is the group of periods
$\Gamma=\text{Span}_{\Z}(1,\alpha,\beta,\gamma)$ (see
Section~\ref{sec-prel}), to the field $\bar{\cK}:= \cK_{\Q}$: we
put
$$Q:= QH_4(M) \otimes_{\cK} \bar{K}\;.$$
In other words, we shall allow any rational powers of the formal
variable $s$ responsible for the symplectic area class. In view of
Proposition \ref{prop-semisimple-algebras-basics}(A) it suffices to
show that $Q$ is semisimple over $\bar{\cK}$.

 We number the facets of $\Delta$ as follows: the facet connecting
$(0, \gamma)$ and $(1-\gamma, \gamma)$ is numbered by $1$, then we
number the rest of the faces moving clockwise from the first one.
The vectors $e_i$ are the following: $e_1 = (0,-1)$, $e_2 =
(-1,-1)$, $e_3=(-1,0)$, $e_4 = (0,1)$, $e_5 = (1,1)$, $e_6 = (1,0)$.
The constants $\eta_1,...,\eta_6$ are given by
$$-\gamma,-1,-\beta,0,\alpha,0\;.$$
We pass to generators of $Q$ by setting $v_i = s^{-\eta_i}qu_i$.
There are nine primitive sets: $\{ 1, 3\}$, $\{ 1, 4\}$, $\{
1,5\}$, $\{ 2, 4\}$, $\{ 2, 5\}$, $\{ 2, 6\}$, $\{ 3, 5\}$, $\{
3,6\}$, $\{ 4,6\}$. One readily checks that the multiplicative
relations are given by $v_iv_{i+2}=v_{i+1}$ and $v_iv_{i+3}=1$,
for all $i=1,...,6$. Here we set $v_i=v_j$ provided
$i=j\;\mod\;6$.  Note that
$$v_2=v_1v_3,v_4=v_1^{-1},v_5=v_1^{-1}v_3^{-1},v_6=v_3^{-1}\;.$$
Put $A=s^{-1/3}v_1, B = s^{-1/3}v_3$. Put $$\epsilon =
\frac{2}{3}-\gamma, \delta= \frac{2}{3}-\beta, \theta= \alpha-
\frac{1}{3}\;.$$ Put
\begin{equation}\label{eq-xyz}
x=s^{\epsilon},y=s^{\delta},z=s^{\theta}\;.
\end{equation}
 Incorporating the
additive relations
$$  -u_1 - u_2 + u_4 + u_5 = 0,\ \ -u_2 - u_3 + u_5 + u_6 = 0$$
we see that $Q=\bar{\cK} [A,B]/I$ where $I$ is the ideal given by
two polynomial equations
\begin{equation}\label{eq-3pts-1}
A^2B^2+xA^2B-B-z=0\;,
\end{equation}
and
\begin{equation}\label{eq-3pts-2}
A^2B^2+yAB^2-A-z=0\;.
\end{equation}
Note that all $u_i$, and hence $A,B$ are invertible. We wish to
show that the ideal $I$ coincides with its radical, which would
yield that $\bar{\cK} [A,B]/ I$ has no nilpotents and is therefore
semi-simple. We shall use the following lemma by Seidenberg, see
\cite{Becker-Weispfenning}, Lemma 8.13: Assume that there exist
polynomials $f_A \in I \cap \bar{\cK}[A]$ and $f_B \in I \cap
\bar{\cK}[B]$ so that $gcd(f_A,f_A')=gcd(f_B,f_B')=1$. Then $I$
coincides with its radical. In view of the symmetry between the
\break $A$- and $B$-variables, it suffices to produce $f_A$ for
all values of $\epsilon, \delta, \theta$.

Multiply equation \eqref{eq-3pts-1} by $(A+y)$ and equation
\eqref{eq-3pts-2} by $A$ and take the difference. We get
\begin{equation} \label{eq-3pts-3}
(A+y)(xA^2-1)B +(A^2-yz)=0\;.
\end{equation}
Multiply this equation by $AB$ and subtract equation
\eqref{eq-3pts-2} multiplied by $(xA^2-1)$. We get
\begin{equation} \label{eq-3pts-4}
A(A^2-yz)B + (A+z)(xA^2-1) =0\;.
\end{equation}

At this point the proof splits into three cases. We shall use
several times the following "non-vanishing test": Let $h$ be a
polynomial of $n$ variables with coefficients in $\C$, and let
$a_1,...,a_n$ be real numbers. Then $t=h(s^{a_1},...,s^{a_n})$ can
be considered as an element of $\bar{\cK}$. Specializing $t$ to
$s=1$ we see that $t\neq 0 \in \bar{\cK}$ provided $h(1,...,1) \neq
0 \in \C$.

Further on, given a polynomial $f$ of the variable $A$ with the
coefficients depending on the parameters $x,y,z$, we denote by
$f'$ the derivative with respect to $A$.

\medskip
\noindent {\sc Case I.} Assume that $y\neq z, xyz \neq 1$.
Multiply equation \eqref{eq-3pts-4} by \break $(A+y)(xA^2-1)$ and
equation \eqref{eq-3pts-3} by $A(A^2-yz)$  and take the
difference. We get that
$$f(A;x,y,z):= (A+y)(A+z)(xA^2-1)^2-A(A^2-yz)^2=0\;.$$
In other words, $f \in I$  as a polynomial of the variable $A$ for
the values of parameters $x,y,z$ given by \eqref{eq-xyz}. We claim
that $gcd(f,f')=1$. Denote by $h(x,y,z) \in \C [x,y,z]$ the
resultant of $f$ and $f'$. It suffices to show that
\begin{equation}
\label{eq-3pts-5} h(s^{\epsilon},s^{\delta},s^{\theta}) \neq 0 \in
\bar{\cK} \end{equation} for all $\epsilon, \delta$ and $\theta$
satisfying assumptions of Case I. Using SINGULAR software (see
\cite{Greuel-Pfister}) we calculate that $h$, interestingly enough,
admits a factorization
$$h(x,y,z)= x^2(y-z)^2(xyz-1)^4h_0(x,y,z),$$
where $h_0$ is some explicit polynomial which we do not write here
as it contains a large number of monomials. Note that $x \neq 0$,
while $y-z$ and $xyz-1$ do not vanish due to our assumptions. Using
SINGULAR again , we calculate that $h_0(1,1,1)=6912\neq 0$ and hence
\eqref{eq-3pts-5} follows. This completes the proof in Case I.

\medskip
\noindent {\sc Case II.} Assume that $y = z, xyz-1 \neq 0$.
Multiply equation \eqref{eq-3pts-3} by $A(A-y)$ and equation
\eqref{eq-3pts-4} by $(xA^2-1)$. Taking into account that $y=z$
and subtracting the obtained equations we get that
$$f(A;x,y):= (A+y)f_0(A;x,y)=0,\; f_0(A;x,y) =
(xA^2-1)^2-A(A-y)^2.$$ Note first that $A=-y$ is not a root of
$f_0$: indeed $f_0(-y;x,y)$ evaluated at $x=s^{\epsilon}$, $y =
s^{\delta}$ and $s=1$ equals to $4 \neq 0$. Thus it suffices to
show that $f_0$ is mutually prime with $f_0'$. Denote by $h(x,y)
\in \C [x,y]$ the resultant of $f_0$ and $f_0'$. It suffices to
show that
\begin{equation}
\label{eq-3pts-6} h(s^{\epsilon},s^{\delta}) \neq 0 \in \bar{\cK}
\end{equation} for all $\epsilon$ and $\delta$ satisfying
assumptions of Case II. Using SINGULAR  we calculate that $h$ admits
a factorization
$$h(x,y,z)= x^2(xy^2-1)^2h_0(x,y),$$
where $h_0 (x,y) = 27x^2y^4+256x^3-192x^2 y-6x y^2-4y^3+27$.
Note that $x \neq 0$, while $xy^2-1$ do not vanish due to our
assumptions. Using SINGULAR, we calculate that $h_0(1,1)=-108\neq
0$ and hence \eqref{eq-3pts-6} follows. This completes the proof
in Case II.

\medskip
\noindent {\sc Case III.} Assume that $xyz=1$.
 Multiply equation
\eqref{eq-3pts-3} by $x^{-1}A$ and equation \eqref{eq-3pts-4} by
$(A+y)$. Taking into account that $yz=x^{-1}$ and subtracting the
obtained equations we get that
$$f(A;x,y):= x(A^2-yz)f_0(A;y,z)=0,$$
$$f_0(A;y,z) = A^2+(y+z-y^2z^2)A+yz.$$ We claim that $f$ and $f'$
are mutually prime. Substituting $A = \pm\sqrt{yz}$ in\-to $f_0$
and specializing to $y=z=1$ we get that $f_0(\pm 1;1,1) = 2\pm 1
\neq 0$. Thus $A^2-yz$ and $f_0$ have no common roots. The
discriminant $d(y,z)$ of $f_0$ equals
$$d(y,z) = (z+y-y^2z^2)^2-4yz\;.$$
Since $d(1,1)=-3\neq 0$ we see that $f_0$ has no multiple roots.
The claim follows. This finishes the proof.
\end{proof}

\medskip
\noindent This finishes the proof of
Theorem~\ref{thm-quantum-hom-sympl-toric-4-mfd-is-semi-simple}.

\section{Quantum homology in the monotone case}
\label{sec-qh-monotone}

We will show now that the semi-simplicity assumption in
Theorem~\ref{thm-qstate-non-monotone} is equivalent in the monotone
case to the one from \cite{EP-qmm}. Assume $M$ is spherically
monotone with $c_1 = \kappa [\omega]$, $\kappa >0$, on $H_2^S (M)$.
In this case one may define the quantum homology ring of $M$ using a
different coefficient ring $\widehat{\Lambda}$, which is a field.
The rings $\widehat{\Lambda}$ and $\Lambda$ are different extensions
of some common subring $\Lambda'$ though none of the rings
$\widehat{\Lambda}$, $\Lambda$ is a subset of another.

The precise definitions are as follows. Denote by $N$ the minimal
Chern number of $M$. Introduce a new variable $w$, and define
$\Lambda'$ as $\Lambda' = \C [w, w^{-1}]$. Identify $\Lambda'$
with a subring of $\Lambda$ through the mapping $w \mapsto
(s^\kappa q)^N$ . The field $\widehat{\Lambda}$, a different ring
extension of $\Lambda'$, is then defined as
$$\widehat{\Lambda} := \bigg\{\ \sum_{\theta \in \Z} z_\theta w^\theta, \ z_\theta \in
\C,\ \sharp\big\{ \theta > c\ |\ z_\theta\neq 0\big\} <\infty,\
\forall c\in\R\ \bigg\}.$$ Define the grading degree of $w$ as
$2N$ and use it to define the grading on $\widehat{\Lambda}$. Note
that $\widehat{\Lambda}$ is a field while $\Lambda'$ is not.

Set $QH'_{ev} (M) := H_{ev} (M; \C) \otimes_\C \Lambda'$. The
monotonicity of $(M,\omega)$ allows to use the Gromov-Witten
invariants in order to define the quantum product on $QH'_{ev}
(M)$ similarly to the case of $QH_* (M)$ in the following way.
Since $(M,\omega)$ is spherically monotone $\bar \pi_2 (M) \cong
\Z$. Denote by $S$ the generator of $\bar \pi_2 (M)$ for which
$\omega (S) >0$.  For $a,b \in H_{ev} (M; \C)$ and $j \in \N$ we
define \break $(a*b)_j \in H_{ev} (M; \C)$ as the unique class
satisfying $(a*b)_j \circ c = GW_{jS} (a,b,c)$ for all $c \in
H_{ev} (M; \C)$, where $\circ$ stands for the ordinary
intersection index in homology and $GW_{jS} (a,b,c)$ is the
Gromov-Witten invariant which, roughly speaking, counts the
pseudo-holomorphic spheres in the homology class $jS$ passing
through cycles representing $a,b,c$. Now for any $a, b\in H_{ev}
(M; \C)$ set
$$a*b = a \cap b + \sum_{j \in \N} (a*b)_j\; w^{-j}.$$
Since $GW_{jS} (a,b,c) = 0$ unless $\deg a + \deg b + \deg c = 4n
- 2Nj$, the sum in the right-hand side is finite and therefore
$a*b\in QH'_{ev} (M)$. By $\C{\hbox{\rm -linearity}}$ extend the
quantum product to the whole $QH'_{ev} (M)$. As a result, for the
same reasons as in the general case, one gets a correctly defined
commutative associative product operation on $QH'_{ev} (M)$
\cite{Liu}, \cite{MS2}, \cite{Ru-Ti}, \cite{Ru-Ti-1}, \cite{Wi}.
It can be easily seen that $QH'_{ev} (M)$ is a subring of $QH_{ev}
(M) = H_{ev} (M; \C) \otimes_\C \Lambda$.

Now $\widehat{QH}_{ev} (M) := H_{ev} (M; \C) \otimes_\C
\widehat{\Lambda}$ and $QH_{ev} (M)$ can be viewed as different
extensions (though not subsets of each other) of the ring
$QH'_{ev} (M)$:
\begin{equation}
\label{eqn-QH-prime} QH_{ev} (M) = QH'_{ev} (M) \otimes_{\Lambda'}
\Lambda,
\end{equation}
\begin{equation}
\label{eqn-QH-prime-prime} \widehat{QH}_{ev} (M) = QH'_{ev} (M)
\otimes_{\Lambda'} \widehat{\Lambda},
\end{equation}
where the tensor product is {\it the tensor product of rings}. Note
that $\widehat{QH}_{ev} (M)$ is a module of finite rank over
$\widehat{\Lambda}$.

Similarly to $QH_* (M)$ the ring $\widehat{QH}_{ev} (M)$ gets
naturally equipped with a grading. This graded ring is precisely
the quantum homology ring used in \cite{EP-qmm}.

\begin{thm}
\label{thm-semisimple-monotone-case} Assume $(M,\omega)$ is a closed
spherically monotone symplectic manifold. Then $QH_{2n} (M)$ is
semi-simple as a $\cK$-algebra if and only if $\widehat{QH}_{ev}
(M)$ is semi-simple as a $\widehat{\Lambda}$-algebra (which is the
semi-simplicity assumption used in \cite{EP-qmm}).

\end{thm}

\begin{proof}
Introduce a new variable $u$ and define a ring $$\cR := \bigg\{\
\sum_{\theta \in \Z} z_\theta u^\theta, \ z_\theta \in \C,\
\sharp\big\{ \theta
> c\ |\ z_\theta\neq 0\big\} <\infty,\ \forall c\in\R\ \bigg\}.$$
The ring $\cR$ is graded with the graded degree of $u$ equal $2$.
Note that $\cR$ is a field. We identify $\widehat{\Lambda}$ with a
subfield of $\cR$ via the mapping $w \mapsto u^N$.

Put $Q:= H_{ev} (M; \C) \otimes_{\C} \cR$ -- this is a graded
$\cR$-algebra that can be viewed as a ring extension of
$\widehat{QH}_{ev} (M)$:
$$Q= \widehat{QH}_{ev} (M)\otimes_{\widehat{\Lambda}} \cR,$$ where the tensor
product is the tensor product of rings. Denote by $Q_{2n}$ the
graded component of degree $2n$ of $Q$. It is a $\C$-algebra.

Recall that the field $\cK$ appearing in the formulation of the
theorem is defined as $\cK_{\omega(S)\Z}$ (see
Section~\ref{sec-prel} above). Note that $\omega(S) = \kappa N$.
Consider an extension $\cK_{\kappa \Z}$ of $\cK = \cK_{\kappa N
\Z}$. Introduce two auxiliary rings
$$R_1 = QH_{2n}(M) \otimes_{\cK} \cK_{\kappa
\Z}\;\;\;\text{and}\;\;\; R_2 = Q_{2n} \otimes_{\C} \cK_{\kappa
\Z}\;.$$ We claim that they are isomorphic as $\cK_{\kappa
\Z}$-algebras. Indeed, notice that every ele\-ment of $R_1$ can be
written as $\sum a_j\otimes q^ms^{\kappa j}$ and every element of
$R_2$ as $\sum a_j\otimes u^ms^{\kappa j}$, where in both cases
the summation conditions are as follows:
$$j \in \Z,\ a_j \in H_{ev}(M,\C),\ \deg a_j +2m =2n,$$
and all $a_j$ vanish for $j$ large enough. The respective
multiplications in $R_1$ and $R_2$ are given by
$$a*_1b= \sum_j (a*b)_{jS}q^{-Nj}s^{-\kappa N j}$$
and
$$a*_2b= \sum_j (a*b)_{jS}u^{-Nj}$$
for all $a,b \in H_{ev}(M,\C)$. Define a map $\Psi: R_2 \to R_1$ by
$$\Psi\Big{(}\sum a_j\otimes u^ms^{\kappa j}\Big{)} = \sum a_j \otimes q^m
s^{\kappa(m+j)}\;,$$ where the summation conditions are as above.
Roughly speaking, $\Psi$ acts by the substitution $u \mapsto
qs^{\kappa}$.  We see that $\Psi$ is a $\cK_{\kappa \Z}$-linear
map. Furthermore, we have $\Psi(a*_2b) = a*_1b$ for all $a,b \in
H_{ev}(M,\C)$. This readily implies that $\Psi$ is a ring
homomorphism. Finally, $\Psi$ is invertible, and the claim
follows.

After these preparations let us pass to the proof of the theorem.
Assume first that $\widehat{QH}_{ev} (M)$ is semi-simple as a
$\widehat{\Lambda}$-algebra. Then, by
Proposition~\ref{prop-semisimple-algebras-basics}(A), $Q$ is a
semi-simple $\cR$-algebra. Decompose $Q$, as a $\cR$-algebra, into a
direct sum of fields $Q_i$. The unity $e_i$ of each such field is an
idempotent and therefore belongs to $Q_{2n}$. It is easy to see that
$Q_{2n}$ decomposes, as a finite-dimensional $\C$-algebra, into a
direct sum of $\C$-algebras $Q_{2n} * e_i$. Each of these algebras
is a field as it lies in $Q_i$. Thus $Q_{2n}$ is a semi-simple
$\C$-algebra. Applying
Proposition~\ref{prop-semisimple-algebras-basics}(A), we get that
since $Q_{2n}$ is $\C$-semi-simple, the ring $R_2$, and hence $R_1$,
is $\cK_{\kappa \Z}$-semi-simple. This in turn yields, again by
Proposition~\ref{prop-semisimple-algebras-basics}(A), the
semi-simplicity of $QH_{2n}(M)$ over $\cK$.

Assume now that $QH_{2n} (M)$ is semi-simple as a $\cK$-algebra.
Then by Proposition~\ref{prop-semisimple-algebras-basics}(A)
$R_1$, and hence $R_2$, is $\cK_{\kappa \Z}$-semi-simple. Applying
again Proposition~\ref{prop-semisimple-algebras-basics}(A) we see
that $Q_{2n}$ is semi-simple as an $\C$-algebra.   Decompose
$Q_{2n}$, as an $\C$-algebra, into a direct sum of fields $Q_i$
with unity elements $e_i$. Note that all the graded components of
$Q$ can be described as $u^i Q_{2n}$, $i\in\Z$. Therefore, as one
can easily see, $Q$ decomposes, as an $\cR$-algebra, into the
direct sum of $\cR$-algebras $\cR*e_i$ each of which is a field.
Thus $Q$ is semi-simple as an $\cR$-algebra and hence,  by
Proposition~\ref{prop-semisimple-algebras-basics}(A),
$\widehat{QH}_{ev} (M)$ is semi-simple as a
$\widehat{\Lambda}$-algebra. This completes the proof.
\end{proof}

\section{Products}

The following theorem shows how to construct new examples of
manifolds with semi-simple quantum homology from the old ones. We
continue working within the class $\cS$ of symplectic manifolds
defined in Section~\ref{sec-fromto} above.

\begin{thm}
\label{thm-product}  Assume that closed symplectic
manifolds
$(M_i^{2n_i},\omega_i)$, $i=1,2$, and $(M,\omega) =
(M_1\times M_2, \omega_1\oplus \omega_2)$ belong to the class
$\cS$ and that at least one of the manifolds $M_i$ satisfies the
condition $H_k (M_i; \C) = 0$ for all odd $k$. Assume also that
each $QH_{2n_i} (M_i)$, $i=1,2$, is semi-simple (as an algebra
over $K_{\Gamma (M_i, \omega_i)})$. Then $QH_{2n_1 + 2n_2} (M)$ is
also semi-simple (as an algebra over $K_{\Gamma (M, \omega)}$).

\end{thm}

\begin{proof} Let us denote for brevity $\Gamma_i
:=\Gamma(M_i,\omega_i)$, $i=1,2$, and $\Gamma:= \Gamma(M,\omega)$.
Note that $\Gamma = \Gamma_1 + \Gamma_2$.

Let $E_i$ be a ring which is also a module over
$\Lambda_{\Gamma_i}$, $i=1,2$. Note that
$$\Lambda_{\Gamma_1+\Gamma_2}=\Lambda_{\Gamma_i} \otimes_{\cK_{\Gamma_i}} \cK_{\Gamma_1+\Gamma_2},\;
i=1,2\;.$$ We define a modified tensor product operation
$\widehat{\otimes}_{\Lambda}$ by
$$E_1 \widehat{\otimes}_{\Lambda} E_2 := \bigg{(}E_1 \otimes_{\cK_{\Gamma_1}}
\cK_{\Gamma_1 +\Gamma_2}\bigg{)} \otimes_{ \Lambda_{\Gamma_1
+\Gamma_2}} \bigg{(}E_2 \otimes_{\cK_{\Gamma_2}} \cK_{\Gamma_1
+\Gamma_2}\bigg{)}.$$ All the tensor products above are
automatically assumed to be tensor products of rings. In simple
words, we extend both rings to $\Lambda_{\Gamma_1
+\Gamma_2}$-modules and consider the usual tensor product over
$\Lambda_{\Gamma_1 +\Gamma_2}$.

Similarly, given algebras $F_i$, $i=1,2$, over $\cK_{\Gamma_i}$
one can define their modified tensor product
$\widehat{\otimes}_{\cK}$, which is an {\it algebra} over
$\cK_\Gamma$, as
$$F_1 \widehat{\otimes}_{\cK} F_2 := \bigg{(}F_1 \otimes_{\cK_{\Gamma_1}}
\cK_{\Gamma_1 +\Gamma_2}\bigg{)} \otimes_{ \cK_{\Gamma_1
+\Gamma_2}} \bigg{(}F_2 \otimes_{\cK_{\Gamma_2}} \cK_{\Gamma_1
+\Gamma_2}\bigg{)},$$ where again all middle tensor products are
tensor products of rings.
Proposition~\ref{prop-semisimple-algebras-basics} implies that if
each $F_i$, $i=1,2$, is a semi-simple algebra over
$\cK_{\Gamma_i}$ then $F_1 \widehat{\otimes}_{\cK} F_2$ is a
semi-simple algebra over $\cK_\Gamma$. As we will show
\begin{equation}
\label{eqn-qhom-hat-cK-product} QH_{2n_1 + 2n_2} (M) = QH_{2n_1}
(M_1) \widehat{\otimes}_{\cK} QH_{2n_2} (M_2),
\end{equation}
which would yield the required semi-simplicity of $QH_{2n_1 +
2n_2} (M)$ as a \break $\cK_\Gamma$-algebra.

Let us prove \eqref{eqn-qhom-hat-cK-product}. Since all the
odd-dimensional homology groups of at least one of the manifolds
$M_1$, $M_2$ are zero, the K\"unneth formula for quantum homology
over the Novikov ring (see e.g. \cite[Exercise 11.1.15]{MS2} for
the statement in the monotone case; the general case in our
algebraic setup can be treated similarly) yields a {\it ring
isomorphism}
$$QH_{ev}(M) =
QH_{ev} (M_1)\widehat{\otimes}_{\Lambda} QH_{ev} (M_2)\;.$$

Note that to each element of $QH_{2n_1+2n_2} (M_1\times M_2)$ of
the form $q^m A\otimes B$ (where $A\otimes B$ denotes the
classical tensor product of two singular homology classes of {\it
even} degree), one can associate in a unique way a pair of
integers $k,l$ with $k+l=m$ so that $q^k A \in QH_{2n_1} (M_1)$
and $q^l B\in QH_{2n_2} (M_2)$. Indeed, the degrees of $A$ and $B$
are even and ${\rm deg}\, q = 2$. Hence we get the required
isomorphism \eqref{eqn-qhom-hat-cK-product}. \end{proof}

\bigskip
The proof above can be easily modified using
Proposition~\ref{prop-direct-summand} to yield the following
result:

\begin{thm}
\label{thm-field-direct-summand-product} Assume that closed
symplectic
manifolds
$(M_i^{2n_i},\omega_i)$, $i=1,2$, and $(M,\omega) =
(M_1\times M_2, \omega_1\oplus \omega_2)$ belong to the class
$\cS$ and that at least one of the manifolds $M_i$ satisfies the
condition $H_k (M_i; \C) = 0$ for all odd $k$. Also assume that
each $QH_{2n_i} (M_i)$, $i=1,2$, contains a field as a direct
summand. Then $QH_{2n_1 + 2n_2} (M)$ also contains a field as a
direct summand.

\end{thm}

\section{McDuff's examples}
\label{sec-mcduff-examples}

In this section we present McDuff's examples of symplectic
manifolds whose quantum homology contains a field as a direct
summand but is not semi-simple \cite{McD-private}. Let $(M^{2n},
\omega)$ be a closed symplectic manifold. Set the base field $\cF
= \C$. Set $p:=[point]\in H_0 (M)$. The manifold $M$ is called
{\it (symplectically) non-uniruled} \cite{McD-Uniruled} if the
genus zero Gromov-Witten invariant $\langle p, a_2,\ldots,
a_k\rangle_{k,\beta}$ of $M$ is zero for any $k\geq 1$ and any
$0\neq \beta\in H_2 (M)$, $a_i\in H_* (M)$. Here the marked points
in the definition of the Gromov-Witten invariant are allowed to
vary freely. In particular, the genus zero 3-point Gromov-Witten
invariants $GW_\beta^{\C} (p,a_2,a_3)$ involved in the definition
of the quantum homology of $M$ over $\C$ (see
Section~\ref{sec-prel}) are zero for all $0\neq \beta\in H_2 (M)$,
$a_2, a_3\in H_* (M)$.

The class of
non-uniruled symplectic manifolds includes the manifolds
satisfying $\left. [\omega]\right|_{\pi_2 (M)} = 0$ or K\"ahler
manifolds with $c_1 <0$ (for instance, complex hypersurfaces in
$\C P^n$ of degree greater than $n+1$).

Let $\tM$ be a one point blow-up of a non-uniruled symplectic
manifold $(M^{2n},\omega)$, $n\geq 2$. Denote by $\tomega$ the
symplectic form on $\tM$, by $E$ the homology class of the
exceptional divisor in $\tM$. We will write $E^i$, $i=1,...,n-1$
for the $i$-th power of $E$ with respect to the classical
cap-product. Introduce the parameter of the blow-up, $\delta :=
[\tomega] (E^{n-1})$.

We wish $M$ and $\tM$ to stay within the class $\cS$ of symplectic
manifolds introduced in Section~\ref{sec-fromto} above. This
holds, for instance, either if $\dim M=4$, $[\omega]$ is rational
and a non-zero rational multiple of $\delta$ lies in the group of
periods $\Gamma(M,\omega)$, or if $\pi_2(M)=0$.

Our goal is to deduce from McDuff's computation of $QH_* (\tM)$
\cite{McD-Uniruled} the following result (an equivalent
formulation was suggested to us by McDuff \cite{McD-private}):

\begin{thm}
\label{thm-mcd-examples-main} If $M$ is non-uniruled, the algebra
 $QH_{2n} (\tM)$ contains a field as a direct summand, but is not
 semi-simple.
\end{thm}

\medskip
\noindent {\bf Convention: In what follows, we use notation
$\dot{+}$ for the direct sum in the category of modules, and
$\oplus$ for the direct sum in the category of algebras.}

\begin{prop}[\cite{McD-Uniruled}, Cor. 2.4]
\label{prop-QH-powers-of-E-tM} Let $0\leq i,j<n$. Then
$$E^i * E^j = E^i\cap E^j = E^{i+j},\ \  {\it if}\ \ i+j<n,$$
$$E^i * E^{n-i} = -p +E\otimes q^{-n+1} s^{-\delta}\;.$$
\end{prop}

\medskip
\noindent Denote $\tcK = \cK_{\Gamma(\tM,\tomega)}$. Set ${\bf 1}
:= [\tM]$, $A:= Eq$, $B:=pq^n$. Let us define the following
$\tcK$-submodules of $QH_{2n} (\tM)$:
$$U := {\rm span}_{\tcK} (A,A^2,\ldots, A^{n-1}),$$
where the powers $A^i$ are taken with respect to the quantum
product,
$$W := QH_{2n} (\tM)\cap \big( \dot{+}_{0<i<2n} H_i (M)\otimes \tLambda \big) .$$
Thus additively
$$QH_{2n} (\tM) = U \dot{+} \tcK B \dot{+} \tcK {\bf 1} \dot{+} W.$$
As for the multiplicative structure of $QH_{2n} (\tM)$, D.McDuff
proved the following

\begin{prop}[\cite{McD-Uniruled}, Prop. 2.5, Cor. 2.6, Lem. 3.3]
\label{prop-QH-tM}
$$U*W =0,$$
$$W*W\subset W\dot{+} \tcK B,$$
$$B* (U\dot{+} W) = 0,$$
$$B*B=0.$$
\end{prop}

\bigskip
 {\sc Proof of Theorem~\ref{thm-mcd-examples-main}.} Set
$V:= U\dot{+} \tcK B \dot{+} \tcK {\bf 1}$. Then one has $QH_{2n}
(\tM) = V\dot{+} W$. By Propositions~\ref{prop-QH-powers-of-E-tM}
and \ref{prop-QH-tM}, $V$ is closed under the quantum product. It
is multiplicatively generated over $\tcK$ by $A$, $B$ and ${\bf
1}$ with the relations:
$$A^n = -B + As^{-\delta},$$
$$A*B=0,$$
see Propositions~\ref{prop-QH-powers-of-E-tM}, \ref{prop-QH-tM}.
Expressing $B$ via $A$ from the first relation and substituting the result into the second one we obtain
the following
isomorphism of $\tcK$-algebras:
$$V\cong \tcK[A]/(f),$$
where $(f)$ is an ideal in $\tcK[A]$
generated by $f(A) := A^2(A^{n-1} - s^{-\delta})$.

Set $g(A):= A^2$, $h(A):= A^{n-1} - s^{-\delta}$, so that $f=gh$.
Decompose $h$ into mutually prime polynomials $P_1, \ldots, P_k\in
\tcK[A]$, where each $P_i$ is irreducible over $\tcK$ (we can
choose $P_i$'s to be irreducible -- and not powers of irreducible
polynomials -- since $h$ has no common factors with its
derivative). By the Chinese remainder theorem, there exists an
isomorphism of $\tcK$-algebras $$I: V=\tcK[A]/(f) \to
\tcK[A]/(g)\oplus \tcK[A]/(P_1)\oplus\ldots\oplus\tcK[A]/(P_k),
$$
where $I$ sends each $u(A) + (f) \in \tcK[A]/(f)$ to the direct
sum of the remainders of the polynomial $u(A)$ modulo $g,
P_1,\ldots, P_k$. Since each $P_i$ is irreducible over $\tcK$, the
factor algebra $\tcK[A]/(P_i)$ is a field. Let $F\subset V$ be a
field which is the preimage of $\tcK[A]/(P_k)$ under the
isomorphism $I$ above. The field $F$ consists of all those
remainders $u\in\tcK[A]$ mod $f$ that are divisible by $g,
P_1,\ldots, P_{k-1}$. In particular, each such $u$ is divisible by
$g(A)=A^2$. In view of Proposition \ref{prop-QH-tM} this yields
\begin{equation}
\label{eq-vsp-W} F*W=0,\; F*B=0\;.
\end{equation}
Put $$Z=
I^{-1}\Big{(}\tcK[A]/(g)\oplus\tcK[A]/(P_1)\oplus\ldots\oplus\tcK[A]/(P_{k-1})\Big{)}\;.$$
Since $V = F \oplus Z$ we have
\begin{equation}\label{eq-vsp-FZ}
F*Z=0,\; Z*Z \subset Z\;.
\end{equation}
Look at the decomposition $$QH_{2n} (\tM) = F\dot{+} Z\dot{+}
W\;.$$ Let us show that  in fact
\begin{equation}
\label{eqn-fzw} QH_{2n} (\tM) = F\oplus( Z\dot{+} W)\;.
\end{equation}
Indeed, by formulas \eqref{eq-vsp-FZ} and \eqref{eq-vsp-W} we have
$F*(Z\dot{+}W)=0$. Thus it suffices to show that $Z\dot{+} W$ is
closed under the quantum multiplication.

Note that $Z*Z \subset Z$ by \eqref{eq-vsp-FZ}. Further, since
$F*B=0$ by \eqref{eq-vsp-W}, the element $B \in V $ necessarily
lies in $Z$. Together with Proposition~\ref{prop-QH-tM} this
implies
\begin{equation}
\label{eqn-WcKB} W*W\subset \tcK B \dot{+} W \subset Z\dot{+} W.
\end{equation}
Finally, $$Z \subset V = U\dot{+} \tcK B \dot{+} \tcK {\bf 1}$$
and hence by Proposition~\ref{prop-QH-tM}
$$Z*W \subset W\;.$$ This yields decomposition \eqref{eqn-fzw}
which tells us that $QH_{2n} (\tM)$ contains a field $F$ as a
direct summand.

At the same time $QH_{2n} (\tM)$ is not semi-simple since it
contains a nilpotent element $B$. This completes the proof. \Qed


\bigskip
{\sc Acknowledgements.} We are grateful to D.~McDuff for
communicating us her unpublished results (\cite{McD-private}, see
Section~\ref{sec-fromto} above) and for encouraging us to include
them into the present paper. We thank M. Borovoi for useful
comments on our proof of Proposition \ref{prop-direct-summand}. We
thank E.~Aljadeff and Y.~Karshon for various useful discussions,
and the anonymous referee for very helpful comments.

\bigskip

\bibliographystyle{alpha}

\bigskip

\setlength{\parindent}{-0.14in}

\begin{tabular}{ll}
Michael Entov & Leonid Polterovich\\
Department of Mathematics&School of Mathematical Sciences\\
Technion - Israel Institute of Technology & Tel Aviv University\\
Haifa 32000, Israel & Tel Aviv 69978, Israel\\
entov@math.technion.ac.il & polterov@post.tau.ac.il\\
\end{tabular}

\end{document}